\documentclass[letterpage, 11pt, notitlepage]{article}

\usepackage[margin=1.0in]{geometry}

\usepackage{epsfig}
\usepackage{amssymb}
\usepackage{amsmath}
\usepackage{amsfonts}
\usepackage{graphicx}
\usepackage{hyperref}
\usepackage{graphicx}
\usepackage{caption}
\usepackage{subcaption}
\usepackage{amsthm,epstopdf,titling,url,array}

\newcommand{\beql}[1]{\begin{equation}\label{#1}}
\newcommand{\eeql}{\end{equation}}
\newcommand{\eqn}[1]{(\ref{#1})}

\theoremstyle{plain}
\newtheorem{thm}{Theorem}
\newtheorem{lem}{Lemma}
\newtheorem{prop}{Proposition}
\newtheorem{cor}{Corollary}
\newtheorem{conj}{Conjecture}

\theoremstyle{definition}

\newtheorem{exmp}{Example}

\theoremstyle{remark}

\begin{document}

\title{\fontsize{25}{25}\selectfont A service system with randomly behaving on-demand agents}
\author{
Lam M. Nguyen\\
Department of Industrial \\ and Systems Engineering\\ 
Lehigh University \\ 
Bethlehem, PA 18015\\
lmn214@lehigh.edu\\
  \and
Alexander L. Stolyar\\
Department of Industrial \\ and Systems Engineering\\ 
Lehigh University \\
Bethlehem, PA 18015\\
stolyar@lehigh.edu\\
}
\maketitle

\begin{center}
\textbf{Abstract}
\end{center}

We consider a service system where agents (or, servers) are invited on-demand. Customers arrive as a Poisson process and join a customer queue. 
Customer service times are i.i.d. exponential. Agents' behavior is random in two respects. First, they can be invited into the system exogenously, and join the agent queue after a random time. Second, with some probability they rejoin the agent queue after a service completion, and otherwise leave the system. 
%The customer and agent queues cannot be non-empty simultaneously -- the head-of-the-line customer and agent are matched immediately and together go to service. 
The objective is to design a real-time adaptive agent invitation scheme that keeps both customer and agent queues/waiting-times small. We study an adaptive scheme, which controls the number of pending agent invitations, based on queue-state feedback.\\

We study the system process fluid limits, in the asymptotic regime where the customer arrival rate goes to infinity. The fluid limit trajectories have complicated behavior -- there are two domains where they follow different ODEs, and a ``reflecting'' boundary. We use the machinery of switched linear systems and common quadratic Lyapunov functions to approach the stability of fluid limits at the desired equilibrium point (with zero queues). We derive sufficient local stability conditions for the fluid limits. We conjecture that, for our model, local stability is in fact sufficient for global stability of fluid limits; the validity of this conjecture is supported by numerical and simulation experiments.
When the local stability conditions do hold, simulations show good overall performance of the scheme. 

\section{Introduction}\label{intro}

We study a service system with exogenously arriving customers, and servers, called {\em agents}, which can be invited to join the system at any time. The system control needs to match the arriving customers with invited agents, with the objective to minimize waiting times of both customers and agents. What makes this problem non-trivial is the fact that there is uncertainty in the agents' behavior. First, invited agents do not arrive into the system immediately; instead they join the system after a random delay. Second, after an agent is done serving a customer, it can either leave the system or return to serve more customers. \\

This model (described in more detail below) is a generalization of that in \cite{stolyar2010pacing,Pang@2014}.
It was originally motivated (see \cite{stolyar2010pacing}) by applications to call/contact centers, where
what we call agents are ``special agents'',
or ``knowledge workers,'' whose time is expensive, so that it is inefficient to have them working fixed shifts, 
with inevitable periods of idle time due to random fluctuations in customer demand. It is much more reasonable 
to invite them on-demand in real time; however, designing an efficient agent invitation strategy is non-trivial 
due to randomness in agent behavior. Besides efficiency (in terms of minimizing customer and agent waiting times),
another highly desirable feature of the invitation scheme is simplicity and robustness.  (For a general discussion of modern call/contact centers and their management, see, e.g. \cite{AAM,liveops} and references therein.)\\

We note that the model we consider is generic and has other applications, or potential applications. One example is telemedicine \cite{WinNT}, in which case ``agents'' are doctors, invited on-demand to serve patients remotely. Another example is crowdsourcing-based customer service \cite{Arise1,Arise2}. Also note that the model has relation to classical assemble-to-order models, where customers are orders and ``invited agents'' are products, which cannot be produced/assembled instantly. The model is also related to ``double-ended queues'' (see e.g. \cite{K66,LGK}) and 
matching systems (see e.g. \cite{Gurvich@2014}); although in such models arrivals of all types into the 
system are typically exogenous, as opposed to being controlled. \\

More specifically, our model is as follows. Customers arrive as a Poisson process and join a customer queue. 
Customer service times are i.i.d. exponential. Agents' behavior is random in two respects. First, they can be invited into the system exogenously, and join the agent queue after a random time. Second, with some probability they rejoin the agent queue after a service completion, and otherwise leave the system. 
(This generalizes the model in \cite{stolyar2010pacing,Pang@2014}, where the agents always leave the system after service completions, thus making our model more realistic in many scenarios.)
The customer and agent queues cannot be non-empty simultaneously -- the head-of-the-line customer and agent are matched immediately and together go to service.
The objective is to design a real-time adaptive agent invitation scheme that keeps both customer and agent queues/waiting-times small. \\

We study a feedback-based adaptive scheme of \cite{stolyar2010pacing,Pang@2014}, which controls the number of pending agent invitations, depending on the customer and/or agent queue lengths and their changes. Due to the fact that our model is more general, the system dynamics is substantially more complicated.\\

The system state can be described by three variables, which are the number of pending invited agents, the difference between agent and customer queues, and the number of customers (or agents) in service. 
For the purposes of analysis, it is more convenient to
consider an alternative, equivalent representation of the system state, which is also described by three variables: the number of pending invited agents, the difference between agent and customer queues, and the total number of customers and agents in the system. \\

We consider the system in the asymptotic regime where the customer arrival rate becomes large while the distributions of an agent response times and a service time are fixed. We show convergence of the fluid-scaled process to the fluid limit (Theorem \ref{thrm1}). The fluid limit trajectories have complicated behavior -- there are two domains where they follow different ODEs, and a ``reflecting'' boundary. This poses big challenges for proving {\em global stability} of the fluid limits, understood as the convergence of their trajectories 
to the equilibrium point, at which the queues are zero.\\

Given that establishing global stability appears to be a very difficult problem,
the focus of this paper and our main results
concern the system {\em local stability} at the equilibrium point, understood as the stability of the 
dynamic system which describes fluid limit trajectories away from the boundary.
We use the machinery of switched linear systems and common quadratic Lyapunov functions \cite{lin@2009,Shorten@2007} to obtain our {\bf main results} (Theorem \ref{thrm2} and \ref{thrm3}),
providing sufficient local stability conditions. 
We conjecture that, for our model, local stability is in fact sufficient for global stability of fluid limits;
 the validity of this conjecture is supported by numerical and simulation experiments.\\

Our simulation experiments also show good overall performance of the feedback scheme  when the local stability conditions do hold. 

\subsection{Organization of the paper}

Section~\ref{sec-notation} contains basic notations, conventions, and abbreviations. Some  background facts on linear systems and switched linear systems are given in Section \ref{necessaryfact}. In Section \ref{model}, we describe the model  in detail. In Section \ref{mainresults} we state the main results of the paper. These results are proved in Sections \ref{fluidscale}, \ref{theorem2proof} and \ref{theorem3proof}. Numerical and simulation experiments are described in Section \ref{numerical}; it also contains our conjectures about global and local stability of fluid limits, supported by these experiments. We conclude  in Section \ref{conclusion}. 

\subsection{Basic notations, conventions and abbreviations} 
\label{sec-notation}

Sets of real and real non-negative numbers are denoted by $\mathbb{R}$ and $\mathbb{R}_{+}$; $\mathbb{R}^d$ and $\mathbb{R}^d_{+}$ are the corresponding vector spaces. The standard Euclidean norm of a vector $x \in \mathbb{R}^n$ is denoted $\|x\|$. For a vector $a$ or matrix $A$, we write their transposes as $a^T$ or $A^T$. We write $x(\cdot)$ to mean the function (or random process) $(x(t), t \geq 0)$. For a real-valued function $x(\cdot): \mathbb{R}_{+} \to \mathbb{R}$, we use either $x^{\prime}(t)$ or $(d/dt)x(t)$ to denote the derivative with respect to $t$, and for $x(\cdot): \mathbb{R}_{+} \to \mathbb{R}^d$, we write $(d/dt)x(t) = (x^{\prime}_1(t),\dots,x^{\prime}_d(t))$. For a real number $x$, let $x^{+} = \max\{x,0\}$ and $x^{-} = - \min\{x,0\}$ and let
\begin{gather*}
\text{sgn}(x) = \begin{cases}
1 \ , \ x > 0 \\
0 \ , \ x = 0 \\
-1 \ , \ x < 0
\end{cases}
\end{gather*}
For $x,y \in \mathbb{R}$, we denote $x \wedge y = \min\{x,y\}$ and $x \vee y = \max\{x,y\}$. Symbol $\Leftrightarrow$ means ``equivalent to''. We write $x^r \to x \in \mathbb{R}^n$ to denote ordinary convergence in $\mathbb{R}^n$. For a finite set of scalar functions $f_n(t)$, $t \geq 0$, $n \in \mathbb{N}$, a point $t$ is called \textit{regular} if for any subset $\mathbb{N}_0 \subseteq \mathbb{N}$, the derivatives 
\begin{gather*}
\frac{d}{dt}\max_{n \in \mathbb{N}_0} f_n(t) \ \text{and} \ \frac{d}{dt}\min_{n \in \mathbb{N}_0} f_n(t)
\end{gather*}
exist. (To be precise, we require that each derivative is proper: both left and right derivatives exist and are equal.) \\

%We use small-$o$ notation for deterministic function: for two real-valued functions $f$ and $g$, we write $f(x) = o(g(x))$ if $\limsup_{x \to \infty} |f(x)/g(x)| = 0$. \\

Abbreviation \textit{u.o.c.} means \textit{uniform on compact sets} convergence of functions, with the argument  determined by the context (usually in $[0,\infty)$); \textit{w.p.1} means \textit{with probability 1}; \textit{i.i.d.} means \textit{independent identically distributed}; RHS means \textit{right hand side}; FSLLN means \textit{functional strong law of large numbers}; \textit{CQLF} means \textit{common quadratic Lyapunov function}. 

\section{Some background facts}\label{necessaryfact}

\subsection{Definitions and results related to switched linear system}

In this paper, we will use some machinery of switched linear systems. Here, we provide some necessary background. Consider a \textit{switched linear system}
\begin{gather}\label{swichedsystem}
\Sigma_S: u^\prime(t) = A(t) u(t) \ , \ A(t) \in \mathcal{A} = \{A_1, \dots, A_m\}
\end{gather}

where $\mathcal{A}$ is a set of matrices in $\mathbb{R}^{n \times n}$, and $t \to A(t)$ is a 
%piecewise constant 
mapping from nonnegative real numbers into $\mathcal{A}$. (Usually, as in \cite{Shorten@2007}, this
mapping is required to be piecewise constant with only finitely many discontinuities in any bounded time-interval.
In our case this additional condition is not important, because our switched system will have a continuous derivative; see equation \eqn{theorem2_system} below.) For $1 \leq i \leq m$, the $i^{th}$ constituent system of the switched linear system (\ref{swichedsystem}) is the \textit{linear time-invariant (LTI) system}
\begin{gather}\label{ltisystem}
\Sigma_{A_i}: u^\prime(t) = A_i u(t). 
\end{gather}

The origin is an \textit{exponentially stable equilibrium} of a switched linear system $\Sigma_s$ if there exist real constants $C > 0$, $a > 0$ such that $\|u(t)\| \leq C e^{-a t} \|u(0)\|$ for $t \geq 0$, for all solutions $u(t)$ of the system (\ref{swichedsystem}) under any $A(t)$ (see \cite{Hespanha@2004,Shorten@2007}). \\
%Murray:1994:MIR:561828

A symmetric square $n \times n$ matrix $M$ with real coefficients is \textit{positive definite} if $z^T M z > 0$ for every non-zero column vector $z \in \mathbb{R}^n$. A symmetric square $n \times n$ matrix $M$ with real coefficients is \textit{negative definite} if $z^T M z < 0$ for every non-zero column vector $z \in \mathbb{R}^n$. A square matrix $A$ is called a \textit{Hurwitz matrix} (or \textit{stable matrix}) if every eigenvalue of $A$ has strictly negative real part (see \cite{pontryagin@1962}). \\

The function $V(u) = u^T P u$ is a \textit{quadratic Lyapunov function} (QLF) for the system $\Sigma_{A} : u^\prime(t) = A u(t)$ if (i) $P$ is symmetric and positive definite, and (ii) $P A + A^T P$ is negative definite. Let $\{A_1,\dots,A_m\}$ be a collection of $n \times n$ Hurwitz matrices, with associated stable LTI systems $\Sigma_{A_1},\dots,\Sigma_{A_m}$. Then the function $V(u) = u^T P u$ is a \textit{common quadratic Lyapunov function} (CQLF) for these systems if $V$ is a QLF for each individual system (see \cite{lin@2009,Shorten@2007}). \\

The following facts will be used in the proof of our results (Theorem \ref{thrm2} and \ref{thrm3}). 

\begin{prop}[\cite{lin@2009,Shorten@2007}]
\label{prop1}
The existence of a CQLF for the LTI systems is sufficient for the exponential stability of a switched linear system. 
\end{prop}

\begin{prop}[\cite{lin@2009,Shorten@2007}]
\label{prop4}
Let $A^{+}$ and $A^{-}$ be Hurwitz matrices in $\mathbb{R}^{n \times n}$, and the difference $A^{+} - A^{-}$ has rank one. Then two systems $u^\prime(t) = A^{+} u(t)$ and $u^\prime(t) = A^{-} u(t)$  have a CQLF if and only if the matrix product $A^{+} A^{-}$ has no negative real eigenvalues. 
\end{prop}

\subsection{Stability of linear systems}

The following facts will also be used in the proof of our results (Theorem \ref{thrm2} and \ref{thrm3}). 

\begin{prop}[\cite{pontryagin@1962}] 
\label{prop2}
Let $L(\lambda) = \det(A - \lambda I) = 0$ be the characteristic equation of matrix $A$:
\begin{gather}
L(\lambda) = a_0 \lambda^3 + a_1 \lambda^2 + a_2 \lambda + a_3 = 0 \ , \ a_0 > 0. 
\end{gather}

Matrix $A$ is Hurwitz if and only if $a_1$, $a_2$, $a_3$ are positive and satisfy $a_1 a_2 > a_0 a_3$. 
\end{prop}

\begin{prop}[\cite{irving@2004}]
\label{prop3}
The general cubic equation has the form
\begin{gather} \label{cubiceq}
a \lambda^3 + b \lambda^2 + c \lambda + d = 0 \ , \ a \neq 0, 
\end{gather}

and discriminant
\begin{gather}\label{discriminant}
\Delta = 18 a b c d - 4 b^3 d + b^2 c^2 - 4 a c^3 - 27 a^2 d^2. 
\end{gather}

If $\Delta > 0$, then the equation has three distinct real roots. 

If $\Delta = 0$, then the equation has a multiple root and all its roots are real. 

If $\Delta < 0$, then the equation has one real root and two nonreal complex conjugate roots.
\end{prop}

\begin{prop}[\cite{shorten@2004}]
\label{prop5}
If $A_1$ is non-singular, the product $A^{-1}_1 A_2$ has no negative eigenvalues if and only if $A_1 + \tau A_2$ is non-singular for all $\tau \geq 0$. 
\end{prop}

\section{Model and algorithm}\label{model}

Our model is a generalization of that considered in \cite{stolyar2010pacing,Pang@2014}.
Customers arrive according to a Poisson process of rate $\Lambda > 0$, and join the customer queue waiting for an available agent and are served in the order of their arrival. There is an infinite pool of potential agents, which can be invited to serve customers. Once being invited, an agent will respond after an independent exponentially distributed random time, with mean $1/\tilde \beta$; it accepts the invitation with probability $a>0$, and otherwise rejects it.
Let $\beta = a \tilde \beta > 0$ be the rate at which an agent accepts the invitation. Agents who accept their invitations join the agent queue, in the order of their arrival. The customer and agent queues cannot be positive simultaneously: the head-of-the-line customer and agents are immediately matched, leave their queues, and together go to service. Each service time is an exponentially distributed random variable with mean $1/\mu$;
after the service completion, the customer leaves the system, while the agent rejoins the agent queue with probability $\alpha\in[0,1)$. Thus, there are two ways in which agents join the queue -- exogenously invited agents accepting invitations and agents already in the system rejoining the queue after service completions.
(The model in \cite{stolyar2010pacing,Pang@2014} is a special case of ours, with $\alpha=0$; in other words, the agents certainly leave the system after service completions, and therefore there is no need to account for agents being in service.) \\

\begin{figure}[h]
  \begin{center}
   \includegraphics[width=0.7\textwidth]{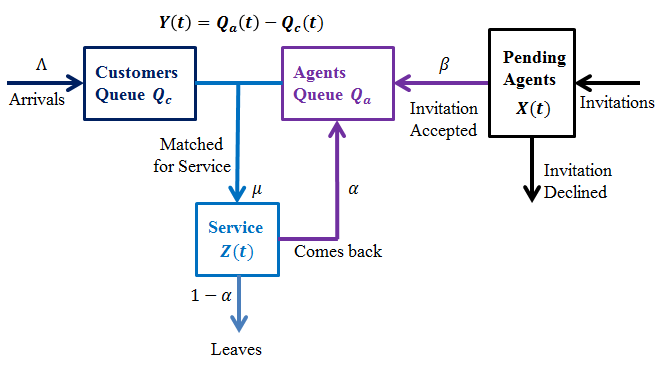}
   \caption{An Agent Invitation System}
   \label{system}
  \end{center}  
\end{figure}

Let $X(t)$ be the number of pending agents that have been invited but have not decided to accept or decline the invitations at time $t$. Let $Q_c(t)$ be the number of customers in the customer queue at time $t$. Let $Q_a(t)$ be the number of agents in the agent queue at time $t$. And we also define $Y(t) = Q_a(t) - Q_c(t)$ as the difference of the agent queue and customer queue at time $t$. Let $Z(t)$ be the number of customers (or agents) in service at time $t$. We assume that the non-idling condition holds, that is, agents do not idle when there are customers waiting in the customer queue, which means that at each time $t$, either the customer queue or the agent queue must be empty. The system state can be described by three variables: $X$: 'the number of pending invited agents'. $Y$: 'the difference between agent and customer queues'. $Z$: 'the number of customers (or agents) in service'. Figure \ref{system} depicts such an agent invitation system. \\

%\cite{stolyar2010pacing,Pang@2014}

The feedback invitation scheme in \cite{stolyar2010pacing}, let us label it as \textit{Scheme A}, is defined as follows. The scheme maintains a ``target'' $X_{target}(t)$ for the number of invited agents $X(t)$. The target $X_{target}(t)$ is changed by $\Delta X_{target}(t) = [-\gamma \Delta Y(t) - \epsilon Y(t) \Delta t]$ at each time $t$ when $Y(t)$ changes by $\Delta Y(t)$ (which can be either $+1$ or $-1$), where $\gamma > 0$ and $\epsilon > 0$ are the algorithm parameters and $\Delta t$ is the time duration from the previous change of $Y$. New agents are invited if and only if $X(t) < X_{target}(t)$, where $X(t)$ is the actual number of invited (pending) agents; therefore, $X(t) \geq X_{target}(t)$ holds at all times. In addition, the target $X_{target}(t)$ is not allowed to go below zero, $X_{target}(t) \geq 0$; i.e. if an update of $X_{target}(t)$ makes it negative, its value is immediately reset to zero. Note that $X_{target}(t)$ is not necessarily an integer. \\

Although the scheme we consider is same as in \cite{stolyar2010pacing}, the model we apply it to is different. Namely, arrivals into the agent queue are not only due to invited agents accepting invitations, but also due to agents returning immediately after the service completions. As a result, the process describing the system evolution contains additional variable $Z$, and is more complicated. \\

To simplify our theoretical analysis, just as in \cite{Pang@2014},
we consider a ``stylized'' version of Scheme A, which has the same basic dynamics, but keeps $X_{target}(t)$ integer and assumes that $X(t) = X_{target}(t)$ at all times; the latter is equivalent to assuming that not only agent invitations can be issued instantly, but they can also be withdrawn at any time. Given these assumptions, when pending agents decline invitations, it has no impact on the system state, because $X(t)$ is immediately ``replenished'' by inviting another agent. Therefore, in the analysis of stylized scheme, the events of declined invitations can be ignored. \\

Formally, the stylized scheme, which we label \textit{Scheme B}, is defined as follows. There are four types of mutually independent, and independent of the past, events that affect the dynamics of $X(t)$, $Y(t)$ and $Z(t)$ in a small time interval $[t, t + dt]$: (i) a customer arrival with probability $\Lambda dt + o(dt)$, (ii) an agent acceptance with probability $\beta X(t) dt + o(dt)$, (iii) an additional event with probability $\epsilon |Y(t)| dt + o(dt)$, and (iv) service completion with probability $\mu Z(t) dt + o(dt)$. 

The changes at these event times are described as follows: 

(i) Upon a customer arrival, if $Y(t) > 0$, $Z(t)$ changes by $\Delta Z(t) = 1$; and if $Y(t) \leq 0$, $Z(t)$ changes by $\Delta Z(t) = 0$. $Y(t)$ changes by $\Delta Y(t) = -1$, and $X(t)$ changes by $\Delta X(t) = \gamma$ (we assume that $\gamma > 0$ is an integer). 

(ii) Upon the acceptance of an invitation, if $Y(t) < 0$, $Z(t)$ changes by $\Delta Z(t) = 1$; and if $Y(t) \geq 0$, $Z(t)$ changes by $\Delta Z(t) = 0$. $Y(t)$ changes by $\Delta Y(t) = 1$, and $X(t)$ changes by $\Delta X(t) = -(\gamma \wedge X(t))$, that is, the change is by $-\gamma$ but $X(t)$ is kept to be nonnegative. 

(iii) Upon the third type of event, if $X(t) \geq 1$, the change $\Delta X(t) = -\text{sgn}(Y(t))$ occurs; and if $X(t) = 0$, the change $\Delta X(t) = 1$ occurs if $Y(t) < 0$ and $\Delta X(t) = 0$ if $Y(t) \geq 0$. 

(iv) Upon the service completion, (a) with probability $\alpha$, if $Y(t) < 0$, the change $\Delta Z(t) = -1 + 1 = 0$ occurs; and if $Y(t) \geq 0$, the change $\Delta Z(t) = -1$ occurs; $Y(t)$ changes by $\Delta Y(t) = 1$, and $\Delta X(t) = -(\gamma \wedge X(t))$. (b) With probability $(1 - \alpha)$, $Z(t)$ changes by $\Delta Z(t) = -1$.  

\section{Main Results}\label{mainresults}

We consider a sequence of systems, indexed by a scaling parameter $r \to \infty$. In the system with index $r$, the arrival rate is $\lambda r$, while the parameters $\alpha$, $\beta$, $\mu$, $\epsilon$, $\gamma$ are constant. The corresponding process is $(X^r, Y^r, Z^r)$, where $X^r = (X^r(t), t \geq 0)$, $Y^r = (Y^r(t), t \geq 0)$ and $Z^r = (Z^r(t), t \geq 0)$. We will center the values of $X^r$, $Y^r$, and $Z^r$ by $\lambda r (1 -\alpha)/\beta$, $0$, and $\lambda r/\mu$, respectively. These values are such that $\beta X^r + \mu \alpha Z^r = \lambda r$, which means that on average the arrival rate of agents into the agent queue matches the rate of customer arrivals. We define fluid-scaled processes with centering
\begin{gather}
\begin{cases} \label{xyz_centered} 
\bar{X}^r = \frac{1}{r}\left(X^r - \frac{\lambda r (1 -\alpha)}{\beta}\right) \\
\bar{Y}^r = \frac{1}{r} Y^r \\
\bar{Z}^r = \frac{1}{r}\left(Z^r - \frac{\lambda r}{\mu}\right). 
\end{cases}
\end{gather}

Let $W$ be the total number of customers and agents in the system. We know that $Y$ is the difference between agent and customer queues (only one of those queues can be positive at any time since we have the non-idling condition) and $Z$ is the number of customers (or agents) in service. From this, $W = |Y| + 2Z$, which is equivalent to $Z = \frac{1}{2}(W - |Y|)$. Instead of using the process $(X, Y, Z)$, we are using a new process $(X, Y, W)$. This process $(X, Y, W)$ is more convenient for the analysis. We have new fluid-scaled processes with centering
\begin{gather}
\begin{cases} \label{xyw_centered} 
\bar{X}^r = \frac{1}{r}\left(X^r - \frac{\lambda r (1 -\alpha)}{\beta}\right) \\
\bar{Y}^r = \frac{1}{r} Y^r \\
\bar{W}^r = \frac{1}{r}\left(W^r - \frac{2 \lambda r}{\mu}\right). 
\end{cases}
\end{gather}

\begin{thm} \label{thrm1}
Consider a sequence of processes $(\bar{X}^r, \bar{Y}^r, \bar{W}^r)$, $r \to \infty$, with deterministic initial states such that $(\bar{X}^r(0), \bar{Y}^r(0), \bar{W}^r(0)) \to (x(0), y(0), w(0))$ for some fixed $(x(0), y(0), w(0)) \in \mathbb{R}^3$, $x(0) \geq -\frac{\lambda (1 -\alpha)}{\beta}$. Then, these processes can be constructed on a common probability space, so that the following holds. W.p.1, from any subsequence of $r$, there exists a further subsequence such that
\begin{gather}
(\bar{X}^r, \bar{Y}^r, \bar{W}^r) \to (x, y, w) \ \ u.o.c. \ \ as \ \ r \to \infty
\end{gather}

where $(x,y,w)$ is a locally Lipschitz trajectory such that at any regular point $t \geq 0$
\begin{gather}
\begin{cases} \label{theorem1_system} 
x^\prime(t) = \begin{cases} 
-\gamma y^\prime(t) - \epsilon y, \ \text{\textit{if}} \ x(t) > -\frac{\lambda (1 -\alpha)}{\beta} \\
[-\gamma y^\prime(t) - \epsilon y] \vee 0, \ \text{\textit{if}} \ x(t) = -\frac{\lambda (1 -\alpha)}{\beta}
\end{cases} \\
y^\prime(t) = \beta x + \frac{1}{2} \alpha \mu (w - |y|) \\
w^\prime(t) = \beta x + \frac{1}{2} (\alpha - 2) \mu (w - |y|). 
\end{cases} 
\end{gather}
\end{thm}

A limit trajectory $(x,y,w)$ specified in Theorem \ref{thrm1} will be called a \textit{fluid limit} starting from $(x(0),y(0),w(0))$. \\

Consider a dynamic system $(x(t),y(t),w(t)) \in \mathbb{R}^3$: 
\begin{gather}\label{theorem2_system}
\begin{cases}  
x^\prime(t) = -\gamma y^\prime(t) - \epsilon y \\
y^\prime(t) = \beta x + \frac{1}{2} \alpha \mu (w - |y|) \\
w^\prime(t) = \beta x + \frac{1}{2} (\alpha - 2) \mu (w - |y|).
\end{cases}
\end{gather}

Note that the RHS of \eqn{theorem2_system} is continuous. \\

This dynamic system describes the dynamics of fluid limit trajectories when the state is away from the boundary $x = -\frac{\lambda (1 -\alpha)}{\beta}$. The {\em non-linear}
system (\ref{theorem2_system}) is a generalization of the linear system, considered in \cite{Pang@2014}. The latter is a special case of (\ref{theorem2_system}) without variable $w$, and with $\alpha = 0$. The system in \cite{Pang@2014} is simply linear, while (\ref{theorem2_system}) has two domains, defined by the sign of $y$. The following results (in Theorem \ref{thrm2} and \ref{thrm3}) provide sufficient exponential stability conditions for the system (\ref{theorem2_system}). 

\begin{thm} \label{thrm2}
\emph{(Sufficient exponential stability condition)}. For any set of positive $\beta$, $\mu$, and $\alpha \in (0,1)$, there exist values of $\gamma > 0$ and $\epsilon > 0$ satisfying the following condition
\begin{gather}\label{stablecond}
\begin{cases}
\frac{\beta \gamma^2}{4} < \epsilon < \frac{\beta \gamma^2}{2} \\
\epsilon > \frac{\beta \gamma^2}{2} - \left(\frac{\alpha \gamma \mu}{2} - \frac{(1 - \alpha) \mu^2}{2 \beta}\right) \\
\gamma > \frac{(1 - \alpha) \mu}{\alpha \beta}. 
\end{cases}
\end{gather}

For the parameters, satisfying this condition, common quadratic Lyapunov function (CQLF) of the system \emph{(\ref{theorem2_system})} exists, and the system \emph{(\ref{theorem2_system})} is exponentially stable. 
\end{thm}

\begin{thm} \label{thrm3}
\emph{(Sufficient exponential stability condition)}. For any set of positive $\beta$, $\mu$, and $\alpha \in (0,1)$, there exist values of $\gamma > 0$ and $\epsilon > 0$ satisfying the following condition
\begin{gather}\label{stablecond2}
\begin{cases}
\epsilon < \frac{\beta \gamma^2}{2} - \frac{\alpha \gamma \mu}{2} \\
\gamma > \frac{\alpha \mu}{\beta}. 
\end{cases}
\end{gather}

For the parameters, satisfying this condition, common quadratic Lyapunov function (CQLF) of the system \emph{(\ref{theorem2_system})} exists, and the system \emph{(\ref{theorem2_system})} is exponentially stable. 
\end{thm}

We say that our fluid-limit system is \textit{globally stable} if every fluid limit trajectory converges to the equilibrium point $(0,0,0)$;
we say that it is \textit{locally stable} if every trajectory of the dynamic system (\ref{theorem2_system}) converges to the equilibrium point $(0,0,0)$. Therefore, the conditions (\ref{stablecond}) and (\ref{stablecond2}) are sufficient for the local stability of our system. We also note that condition (\ref{stablecond2}) is more robust and is easier to achieve in practice. Indeed, for any given $\epsilon > 0$, it holds for all sufficiently large $\gamma$; how large, can be determined if some estimates of other parameters are available. 

\section{Fluid scale analysis and proof of Theorem \ref{thrm1}}\label{fluidscale}

In order to prove Theorem \ref{thrm1}, it suffices to show that w.p.1 from any subsequence of $r$, we can choose a further subsequence, along which a u.o.c. convergence to a fluid limit holds. \\

Given the initial state $(X^r(0),Y^r(0),W^r(0))$, we construct the process $(X^r, Y^r, W^r)$, for all $r$, on the same probability space via a common set of independent Poisson process as follows:
\begin{gather} 
\label{pro_x} X^r(t) = G^r(t) + \left( -\min_{0 \leq s \leq t} G^r(s)\right) \vee 0, \\
G^r(t) = X^r(0) + \gamma N_1(\lambda rt) - \gamma N_2 \left(\beta \int_0^t X^r(s) ds\right) - \gamma N_4 \left(\alpha \mu \int_0^t \frac{1}{2} (W^r(s) - |Y^r(s)|) ds\right) + \nonumber \\
\label{pro_g} + N_5 \left(\epsilon \int_0^t (Y^r(s))^{-}ds\right) - N_6 \left(\epsilon \int_0^t (Y^r(s))^{+}ds\right),
\end{gather} \begin{gather}
\label{pro_y} Y^r(t) = Y^r(0) + N_2 \left(\beta \int_0^t X^r(s)ds\right) + N_4 \left(\alpha \mu \int_0^t \frac{1}{2} (W^r(s) - |Y^r(s)|) ds\right) - N_1(\lambda rt), \\
W^r(t) = W^r(0) + N_1 (\lambda r t) + N_2 \left(\int_0^t \beta X^r(s) ds\right) - N_3 \left(\int_0^t 2 (1 - \alpha) \mu \frac{1}{2} (W^r(s) - |Y^r(s)|) ds\right) - \nonumber \\
\label{pro_w} - N_4 \left(\int_0^t \alpha \mu \frac{1}{2} (W^r(s) - |Y^r(s)|) ds\right),
\end{gather}

and $N_i(\cdot)$, $i = 1, \dots, 6$ are mutually independent unit-rate Poisson processes \cite{Pang@2007}. $N_1$ is the process which drives customer arrivals. $N_2$ is the process which drives the acceptance of invitations. $N_3$ is the process which drives the service completions, with agents leaving the system. $N_4$ is the process which drives the service completions, with agents coming back. $N_5$ and $N_6$ are the processes which drive the third type of event. W.p.1, for any $r$, relations (\ref{pro_x})-(\ref{pro_w}) uniquely define the realization of $(X^r,Y^r,W^r)$ via the realizations of the driving processes $N_i(\cdot)$. Relation (\ref{pro_x}), the ``reflection'' at zero,  corresponds to the property that $X^r(t)$ cannot become negative. \\

The functional strong law of large numbers (FSLLN) 
%\cite{opac-b1085509} 
holds for each Poisson process $N_i$: 
\begin{gather} \label{fslln}
\frac{N_i(rt)}{r} \to t \ , \ r \to \infty \ , \ \text{u.o.c.}, \ \text{w.p.1}.
\end{gather}

We consider the sequence of associated fluid-scaled processes $(\bar{X}^r,\bar{Y}^r,\bar{W}^r)$ as defined in (\ref{xyw_centered}). (Note that these processes are centered.) Let a constant $m > \|(x(0),y(0),w(0)\|$ be fixed. For each $r$, on the same probability space as $(\bar{X}^r,\bar{Y}^r,\bar{W}^r)$, let us define a modified fluid-scaled process $(\bar{X}^r_m,\bar{Y}^r_m,\bar{W}^r_m)$. Let $(\bar{X}^r_m,\bar{Y}^r_m,\bar{W}^r_m)$ start from the same initial state as $(\bar{X}^r,\bar{Y}^r,\bar{W}^r)$  , i.e., $(\bar{X}^r_m(0),\bar{Y}^r_m(0),\bar{W}^r_m(0)) = (\bar{X}^r(0),\bar{Y}^r(0),\bar{W}^r(0))$. The modified process  $(\bar{X}^r_m,\bar{Y}^r_m,\bar{W}^r_m)$ follows the same path as $(\bar{X}^r,\bar{Y}^r,\bar{W}^r)$ until the first time that $\|(\bar{X}^r(t),\bar{Y}^r(t),\bar{W}^r(t))\| \geq m$. Denote this time by $\tau^r_m$. We then freeze the process $(\bar{X}^r_m,\bar{Y}^r_m,\bar{W}^r_m)$ at the value $(\bar{X}^r(\tau^r_m),\bar{Y}^r(\tau^r_m),\bar{W}^r(\tau^r_m))$, i.e. $(\bar{X}^r_m(t),\bar{Y}^r_m(t),\bar{W}^r_m(t)) = (\bar{X}^r(\tau^r_m),\bar{Y}^r(\tau^r_m),\bar{W}^r(\tau^r_m))$ for all $t \geq \tau^r_m$. 

\begin{lem}
Fix $(x(0),y(0),w(0))$ and a finite constant $m > \|(x(0),y(0),w(0))\|$. Then, w.p.1  for any subsequence of $r$, there exists a further subsequence, along which $(\bar{X}^r_m, \bar{Y}^r_m, \bar{W}^r_m)$ converges u.o.c. to a Lipschitz continuous trajectory $(x_m,y_m,w_m)$, which satisfies properties \emph{(\ref{theorem1_system})} at any regular time $t \geq 0$ such that $\|(x_m(t),y_m(t),w_m(t))\| < m$. 
\end{lem}

\textit{Proof}. For the modified fluid-scaled processes $(\bar{X}^r_m, \bar{Y}^r_m, \bar{W}^r_m)$, we define the associated counting processes for upward and downward jumps. For $t \leq \tau^r_m$, 
\begin{gather} 
\bar{X}^{r \uparrow}_m(t) = r^{-1} \gamma N_1(\lambda rt) + r^{-1} N_5 \left(\epsilon r \int_0^t (\bar{Y}^r_m(s))^{-}ds\right), \\
\bar{X}^{r \downarrow}_m(t) = r^{-1} \gamma N_2 \left(\beta r \int_0^t \left[\bar{X}^r_m(s) + \frac{\lambda (1 -\alpha)}{\beta}\right] ds\right) + \nonumber \\ 
+ r^{-1} \gamma N_4 \left(\frac{1}{2} \alpha \mu r \int_0^t \left[\bar{W}^r_m(s) + \frac{2 \lambda}{\mu} - |\bar{Y}^r_m(s)|\right] ds\right) + r^{-1} N_6 \left(\epsilon r \int_0^t (\bar{Y}^r_m(s))^{+}ds\right), \\
\bar{Y}^{r \uparrow}_m(t) = r^{-1} N_2 \left(\beta r \int_0^t \left[\bar{X}^r_m(s) + \frac{\lambda (1 -\alpha)}{\beta}\right] ds\right) + r^{-1} N_4 \left(\frac{1}{2} \alpha \mu r \int_0^t \left[\bar{W}^r_m(s) + \frac{2 \lambda}{\mu} - |\bar{Y}^r_m(s)|\right] ds\right), \\
\bar{Y}^{r \downarrow}_m(t) = r^{-1} N_1(\lambda rt),
\end{gather} \begin{gather}
\bar{W}^{r \uparrow}_m(t) = r^{-1} N_1 (\lambda r t) + r^{-1} N_2 \left(\beta r \int_0^t \left[\bar{X}^r_m(s) + \frac{\lambda (1 -\alpha)}{\beta}\right] ds\right),  \\ 
\bar{W}^{r \downarrow}_m(t) = r^{-1} N_3 \left((1 - \alpha) \mu r \int_0^t \left[\bar{W}^r_m(s) + \frac{2 \lambda}{\mu} - |\bar{Y}^r_m(s)|\right] ds\right) + \nonumber \\
+ r^{-1} N_4 \left(\frac{1}{2} \alpha \mu r \int_0^t \left[\bar{W}^r_m(s) + \frac{2 \lambda}{\mu} - |\bar{Y}^r_m(s)|\right] ds\right),
\end{gather}

and for $t > \tau^r_m$, all these counting processes are frozen at their values at time $\tau^r_m$, that is, 
\begin{gather}
\begin{cases}
\bar{X}^{r \uparrow}_m(t) = \bar{X}^{r \uparrow}_m(\tau^r_m) \ , \ \bar{X}^{r \downarrow}_m(t) = \bar{X}^{r \downarrow}_m(\tau^r_m) \ , \\ 
\bar{Y}^{r \uparrow}_m(t) = \bar{Y}^{r \uparrow}_m(\tau^r_m) \ , \ \bar{Y}^{r \downarrow}_m(t) = \bar{Y}^{r \downarrow}_m(\tau^r_m) \ , \\
\bar{W}^{r \uparrow}_m(t) = \bar{W}^{r \uparrow}_m(\tau^r_m) \ , \ \bar{W}^{r \downarrow}_m(t) = \bar{W}^{r \downarrow}_m(\tau^r_m). 
\end{cases}
\end{gather}

Using the relations (\ref{pro_x})-(\ref{pro_w}) and the fact that for $0 \leq t \leq \tau^r_m$ the original process $(\bar{X}^r,\bar{Y}^r,\bar{W}^r)$ and the modified process $(\bar{X}^r_m,\bar{Y}^r_m,\bar{W}^r_m)$ coincide, we have for all $t \geq 0$, 
\begin{gather} 
\bar{X}^r_m(t) = \bar{G}^r_m(t) + \left(-\lambda (1 -\alpha)/ \beta - \min_{0 \leq s \leq t} \bar{G}^r_m(s)\right) \vee 0, \\
\bar{G}^r_m(t) = \bar{X}^r(0) + \bar{X}^{r \uparrow}_m(t) - \bar{X}^{r \downarrow}_m(t), \\
\bar{Y}^r_m(t) = \bar{Y}^r(0) + \bar{Y}^{r \uparrow}_m(t) - \bar{Y}^{r \downarrow}_m(t), \\
\bar{W}^r_m(t) = \bar{W}^r(0) + \bar{W}^{r \uparrow}_m(t) - \bar{W}^{r \downarrow}_m(t). 
\end{gather}

The counting processes $\bar{X}^{r \uparrow}_m$, $\bar{X}^{r \downarrow}_m$, $\bar{Y}^{r \uparrow}_m$, $\bar{Y}^{r \downarrow}_m$, $\bar{W}^{r \uparrow}_m$, $\bar{W}^{r \downarrow}_m$ are non-decreasing. Using the Functional Strong Law of Large Number (FSLLN) (\ref{fslln}) and the fact that the processes $\bar{X}^r_m$, $\bar{Y}^r_m$, and $\bar{W}^r_m$ are uniformly bounded by construction, we see that w.p.1. for any subsequence of $r$, there exists a further subsequence along which the set of trajectories $(\bar{X}^{r \uparrow}_m, \bar{X}^{r \downarrow}_m, \bar{Y}^{r \uparrow}_m, \bar{Y}^{r \downarrow}_m, \bar{W}^{r \uparrow}_m, \bar{W}^{r \downarrow}_m)$ converges u.o.c. to a set of non-decreasing  Lipschitz continuous functions $(x^{\uparrow}_m, x^{\downarrow}_m, y^{\uparrow}_m, y^{\downarrow}_m, w^{\uparrow}_m, w^{\downarrow}_m)$. But then the u.o.c. convergence of $(\bar{X}^r_m, \bar{Y}^r_m, \bar{W}^r_m, \bar{G}^r_m)$ to a set of Lipschitz continuous functions $(x_m, y_m, w_m, g_m)$ holds, where
\begin{gather} 
x_m(t) = g_m(t) + \left(-\lambda (1 -\alpha)/ \beta - \min_{0 \leq s \leq t} g_m(s)\right) \vee 0, \\
g_m(t) = x(0) + x^{\uparrow}_m(t) - x^{\downarrow}_m(t), \\
y_m(t) = y(0) + y^{\uparrow}_m(t) - y^{\downarrow}_m(t), \\
w_m(t) = w(0) + w^{\uparrow}_m(t) - w^{\downarrow}_m(t),
\end{gather}

and the following holds for $t$ before fluid trajectory hits $\|(x_m(t),y_m(t),w_m(t))\| = m$
\begin{gather} 
x^{\uparrow}_m(t) = \gamma \lambda t + \epsilon \int_0^t y^{-}_m(s) ds, \\
x^{\downarrow}_m(t) = \gamma \beta \int_0^t \left(x_m(s) + \frac{\lambda (1 -\alpha)}{\beta}\right) ds + \frac{1}{2} \gamma \alpha \mu \int_0^t \left(w_m(s) + \frac{2 \lambda}{\mu} - |y_m(s)|\right) ds + \epsilon \int_0^t y^{+}_m(s) ds,  \\
y^{\uparrow}_m(t) = \beta \int_0^t \left(x_m(s) + \frac{\lambda (1 -\alpha)}{\beta}\right) ds + \frac{1}{2} \alpha \mu \int_0^t \left(w_m(s) + \frac{2 \lambda}{\mu} - |y_m(s)|\right) ds, \\
y^{\downarrow}_m(t) = \lambda t,
\end{gather} \begin{gather}
w^{\uparrow}_m(t) = \lambda t + \beta \int_0^t \left(x_m(s) + \frac{\lambda (1 -\alpha)}{\beta}\right) ds,  \\ 
w^{\downarrow}_m(t) = (1 - \alpha)\mu \int_0^t \left(w_m(s) + \frac{2 \lambda}{\mu} - |y_m(s)|\right) ds + \frac{1}{2} \alpha \mu \int_0^t \left(w_m(s) + \frac{2 \lambda}{\mu} - |y_m(s)|\right) ds.
\end{gather}

It is easy to verify that, for $t$ before fluid trajectory hits $\|(x_m(t),y_m(t),w_m(t))\| = m$
\begin{gather}
\begin{cases}  
x^{\prime}_m(t) = \begin{cases} 
-\gamma \beta x_m - \frac{1}{2} \gamma \alpha \mu w_m + \frac{1}{2} \gamma \alpha \mu |y_m| - \epsilon y_m, \ \text{if} \ x_m(t) > -\frac{\lambda (1 -\alpha)}{\beta} \\
[-\gamma \beta x_m - \frac{1}{2} \gamma \alpha \mu w_m + \frac{1}{2} \gamma \alpha \mu |y_m| - \epsilon y_m] \vee 0, \ \text{if} \ x_m(t) = -\frac{\lambda (1 -\alpha)}{\beta}
\end{cases} \\
y{^\prime}_m(t) = \beta x_m + \frac{1}{2} \alpha \mu (w_m - |y_m|) \\
w{^\prime}_m(t) = \beta x_m + \frac{1}{2} (\alpha - 2) \mu (w_m - |y_m|)
\end{cases}
\end{gather}

which is equivalent to
\begin{gather}
\begin{cases}  
x{^\prime}_m(t) = \begin{cases} 
-\gamma y{^\prime}_m(t) - \epsilon y_m, \ \text{if} \ x_m(t) > -\frac{\lambda (1 -\alpha)}{\beta} \\
[-\gamma y^{\prime}_m(t) - \epsilon y_m] \vee 0, \ \text{if} \ x_m(t) = -\frac{\lambda (1 -\alpha)}{\beta}
\end{cases} \\
y{^\prime}_m(t) = \beta x_m + \frac{1}{2} \alpha \mu (w_m - |y_m|) \\
w{^\prime}_m(t) = \beta x_m + \frac{1}{2} (\alpha - 2) \mu (w_m - |y_m|). 
\end{cases}
\end{gather}

This means properties (\ref{theorem1_system}) hold for the trajectory $(x_m,y_m,w_m)$. This completes the proof. $\Box$ \\

\textit{Conclusion of the proof of Theorem \ref{thrm1}}. It is obvious that inequality $\frac{d}{dt}\|(x_m(t),y_m(t),w_m(t))\| \leq C \|(x_m(t),y_m(t),w_m(t))\|$ holds for any $m$, and some common $C > 0$. From Gronwall's inequality \cite{opac-b1080363}, we have $\|(x_m(t),y_m(t),w_m(t))\| \leq \|(x(0),y(0),w(0))\| e^{C t}$ for $t \geq 0$. For a given $(x(0),y(0),w(0))$, let us fix $T_l > 0$ and choose $m_l > \|(x(0),y(0),w(0)\| e^{C T_l}$. For this $T_l > 0$, there exists a subsequence $r^{l}$, along which $(\bar{X}^r, \bar{Y}^r, \bar{W}^r)$ converges uniformly to $(x_{m_l},y_{m_l},w_{m_l})$, which satisfies properties (\ref{theorem1_system}), at any $t \in [0,T_l]$. The limit trajectory $(x_{m_l},y_{m_l},w_{m_l})$ does not hit $m_l$ in $[0,T_l]$. Subsequence $r^{l} = \{r^{l}_1, r^{l}_2, \dots\}$ is such that, w.p.1, for all sufficiently large $r$ along the subsequence $r^{l}$, $(\bar{X}^r(t),\bar{Y}^r(t),\bar{W}^r(t)) = (\bar{X}^r_{m_l}(t),\bar{Y}^r_{m_l}(t),\bar{W}^r_{m_l}(t))$ at any $t \in [0,T_l]$. We consider a sequence $T_1$, $T_2$, $\dots$, $\to \infty$. We construct a subsequence $r^{*}$ by using Cantor's diagonal process \cite{opac-b1098274} from subsequences $r^{1}$, $r^{2}$, $\dots$ ($r^{1} \supseteq r^{2} \supseteq \dots$) corresponding to $T_1$, $T_2$, $\dots$, respectively (i.e. $r^{*}_1 = r^{1}_1$, $r^{*}_2 = r^{2}_2$, $\dots$). Clearly, for this subsequence $r^{*}$, w.p.1, $(\bar{X}^r, \bar{Y}^r, \bar{W}^r)$ converges u.o.c. to $(x,y,w)$, which satisfies properties (\ref{theorem1_system}), at any regular point $t \in [0,\infty)$. $\Box$ 

\section{Proof of Theorem \ref{thrm2}}\label{theorem2proof}

We use the machinery of switched linear systems and common quadratic Lyapunov functions (CQLF) to approach the stability of fluid limits, i.e. their convergence to the unique equilibrium point $(0,0,0)$ \cite{lin@2009,Shorten@2007}. \\

System (\ref{theorem2_system}) is a switched linear system with $m = 2$. Namely, for $y \geq 0$,  
\begin{gather}\label{system_aplus}
\begin{cases}  
x^\prime(t) = 
\left(-\gamma \beta \right) x + \left(\frac{1}{2} \gamma \alpha \mu  - \epsilon \right) y + \left(-\frac{1}{2} \gamma \alpha \mu \right) w
 \\
y^\prime(t) = \left(\beta \right) x + \left(-\frac{1}{2} \alpha \mu \right) y + \left(\frac{1}{2} \alpha \mu \right) w \\
w^\prime(t) = \left(\beta \right) x + \left(-\frac{1}{2} (\alpha - 2) \mu \right) y + \left(\frac{1}{2} (\alpha - 2) \mu \right) w
\end{cases}
\end{gather}

and for $y < 0$,  
\begin{gather}\label{system_aminus}
\begin{cases}  
x^\prime(t) = 
\left(-\gamma \beta \right) x + \left(-\frac{1}{2} \gamma \alpha \mu  - \epsilon \right) y + \left(-\frac{1}{2} \gamma \alpha \mu \right) w
 \\
y^\prime(t) = \left(\beta \right) x + \left(\frac{1}{2} \alpha \mu \right) y + \left(\frac{1}{2} \alpha \mu \right) w \\
w^\prime(t) = \left(\beta \right) x + \left(\frac{1}{2} (\alpha - 2) \mu \right) y + \left(\frac{1}{2} (\alpha - 2) \mu \right) w. 
\end{cases}
\end{gather}

We can rewrite the systems above as two linear time-invariant systems $u^\prime(t) = A^{+} u(t)$ and $\ u^\prime(t) = A^{-} u(t)$, where $u(t) = (x(t),y(t),w(t))^T$ and
\begin{gather}\label{aplus}
A^{+} = \left( \begin{array}{ccc}
-\gamma \beta & \frac{1}{2} \gamma \alpha \mu  - \epsilon & -\frac{1}{2} \gamma \alpha \mu \\
\beta & -\frac{1}{2} \alpha \mu & \frac{1}{2} \alpha \mu \\
\beta & -\frac{1}{2} (\alpha - 2) \mu & \frac{1}{2} (\alpha - 2) \mu \end{array} \right)
\end{gather} 

and
\begin{gather}\label{aminus}
A^{-} = \left( \begin{array}{ccc}
-\gamma \beta & -\frac{1}{2} \gamma \alpha \mu  - \epsilon & -\frac{1}{2} \gamma \alpha \mu \\
\beta & \frac{1}{2} \alpha \mu & \frac{1}{2} \alpha \mu \\
\beta & \frac{1}{2} (\alpha - 2) \mu & \frac{1}{2} (\alpha - 2) \mu \end{array} \right). 
\end{gather}

\begin{lem} \label{lemmaaplus}
Matrix $A^{+}$ in \emph{(\ref{aplus})} is Hurwitz for all positive $\beta$, $\gamma$, $\mu$, $\epsilon$ and $\alpha \in (0,1)$. 
\end{lem}

\textit{Proof}. The characteristic equation of $A^{+}$ is $\det(A^{+} - \lambda I) = 0$, which is equivalent to
\begin{gather}
\lambda^3 + (\beta \gamma + \mu) \lambda^2 + (\beta \epsilon + \beta \gamma \mu) \lambda + \beta \epsilon \mu = 0. 
\end{gather}

By Proposition \ref{prop2}, it suffices to verify that
\begin{gather}\label{aplus1}
\beta \gamma + \mu > 0 \ , \ \beta \epsilon + \beta \gamma \mu > 0 \ , \ \beta \epsilon \mu > 0, 
\end{gather} 

and
\begin{gather}\label{aplus2}
(\beta \gamma + \mu)(\beta \epsilon + \beta \gamma \mu) - \beta \epsilon \mu = \beta^2 \gamma^2 \mu + \beta^2 \gamma \epsilon + \beta \gamma \mu^2 > 0. 
\end{gather}

Conditions (\ref{aplus1}) and (\ref{aplus2}) are obviously true. $\Box$  

\begin{lem} \label{lemmaaminus}
Matrix $A^{-}$ in \emph{(\ref{aminus})} is Hurwitz for positive $\beta$, $\gamma$, $\mu$, $\epsilon$, and $\alpha \in (0,1)$, satisfying
\begin{gather}\label{aminushurwitz}
\left(\frac{\beta \gamma}{\mu} + (1 - \alpha)\right)\left(\frac{\gamma \mu}{\epsilon} + 1\right) > 1. 
\end{gather}
\end{lem}
\textit{Proof}. The characteristic equation of $A^{-}$ is $\det(A^{-} - \lambda I) = 0$, which is equivalent to
\begin{gather}
\lambda^3 + (\beta \gamma + \mu(1 - \alpha)) \lambda^2 + (\beta \epsilon + \beta \gamma \mu) \lambda + \beta \epsilon \mu = 0. 
\end{gather}

By Proposition \ref{prop2}, it suffices to verify that
\begin{gather}\label{aminus1}
\beta \gamma + \mu(1 - \alpha) > 0 \ , \ \beta \epsilon + \beta \gamma \mu > 0 \ , \ \beta \epsilon \mu > 0,
\end{gather} 

and 
\begin{gather*}
(\beta \gamma + \mu(1 - \alpha))(\beta \epsilon + \beta \gamma \mu) - \beta \epsilon \mu > 0 \ \text{which is equivalent to} \ \left(\frac{\beta \gamma}{\mu} + (1 - \alpha)\right)\left(\frac{\gamma \mu}{\epsilon} + 1\right) > 1. 
\end{gather*}

Conditions (\ref{aminus1}) are obviously true. $\Box$ 

\begin{lem}\label{aminuscond}
For $\beta > 0$, $\mu > 0$ and $\alpha \in (0,1)$, there exists a pair of $\gamma > 0$ and $\epsilon > 0$ satisfying condition
\begin{gather}\label{lemma5}
\begin{cases}
\frac{\beta \gamma^2}{4} < \epsilon < \frac{\beta \gamma^2}{2} \\
\epsilon > \frac{\beta \gamma^2}{2} - \left(\frac{\alpha \gamma \mu}{2} - \frac{(1 - \alpha) \mu^2}{2 \beta}\right) \\
\gamma > \frac{(1 - \alpha) \mu}{\alpha \beta}. 
\end{cases}
\end{gather}

Moreover, condition \emph{(\ref{lemma5})} implies matrix $A^{-}$ being Hurwitz. 
\end{lem}

\textit{Proof}. For $\beta > 0$, $\mu > 0$ and $\alpha \in (0,1)$, we have $\frac{(1 - \alpha) \mu}{\alpha \beta} > 0$. Hence, we can always find a value of $\gamma > 0$ satisfying the third condition of (\ref{lemma5}). And from the third condition of (\ref{lemma5}), we have 
\begin{gather}
\frac{\alpha \gamma \mu}{2} - \frac{(1 - \alpha) \mu^2}{2 \beta} > 0. 
\end{gather}

Hence, we can always find a value of $\epsilon > 0$ satisfying
\begin{gather}
\begin{cases}
\frac{\beta \gamma^2}{4} < \epsilon < \frac{\beta \gamma^2}{2} \\
\epsilon > \frac{\beta \gamma^2}{2} - \left(\frac{\alpha \gamma \mu}{2} - \frac{(1 - \alpha) \mu^2}{2 \beta}\right). 
\end{cases}
\end{gather}

As shown in the proof of Lemma \ref{lemmaaminus}, matrix $A^{-}$ when
\begin{gather*}
(\beta \gamma + \mu(1 - \alpha))(\beta \epsilon + \beta \gamma \mu) - \beta \epsilon \mu = \beta^2 \gamma \epsilon + \beta \epsilon \mu(1 - \alpha) + \beta^2 \gamma^2 \mu + \beta \gamma \mu^2 (1 - \alpha) - \beta \epsilon \mu   > 0. 
\end{gather*}

For positive $\beta$, $\gamma$, $\mu$, $\epsilon$, and $\alpha \in (0,1)$, the condition $\beta^2 \gamma^2 \mu - \beta \epsilon \mu  > 0$ or, equivalently, $\epsilon < \gamma^2 \beta$ implies (\ref{aminushurwitz}). It means that, if $\epsilon < \gamma^2 \beta$, then $A^{-}$ is Hurwitz. But, the condition (\ref{lemma5}) implies $\epsilon < \gamma^2 \beta$. $\Box$ \\

\textit{Conclusion of the proof of Theorem \ref{thrm2}}. The characteristic equation of $A^{+} A^{-}$ is
\begin{gather} \label{eqaplusaminus}
\lambda^3 - (\mu^2 - \alpha \mu^2 + \beta^2 \gamma^2 - 2 \beta \epsilon - \alpha \beta \gamma \mu) \lambda^2 + (\beta^2 \epsilon^2 + \beta^2 \gamma^2 \mu^2 - 2 \beta \epsilon \mu^2 + \alpha \beta \epsilon \mu^2) \lambda - \beta^2 \epsilon^2 \mu^2 = 0. 
\end{gather}

(Expression (\ref{eqaplusaminus}) is obtained with the help of MATLAB symbolic calculation.) By Proposition \ref{prop3}, if $\Delta < 0$, then the equation has one real root and two nonreal complex conjugate roots. It is well known that the determinant of a square matrix $A^{+} A^{-}$ is the product of its eigenvalues. We have $\det(A^{+} A^{-}) = \lambda_1 \lambda_2 \lambda_3 = \beta^2 \epsilon^2 \mu^2 > 0$. Therefore, one of the roots must be a real positive. We see that it will suffice to show that $\Delta < 0$ to demonstrate $A^{+} A^{-}$ could have no negative real eigenvalues. From (\ref{eqaplusaminus}), we have 
\begin{gather}
\begin{cases}
a = 1 \\
b = - (\mu^2 - \alpha \mu^2 + \beta^2 \gamma^2 - 2 \beta \epsilon - \alpha \beta \gamma \mu) \\
c = \beta^2 \epsilon^2 + \beta^2 \gamma^2 \mu^2 - 2 \beta \epsilon \mu^2 + \alpha \beta \epsilon \mu^2 \\
d = - \beta^2 \epsilon^2 \mu^2. 
\end{cases}
\end{gather}

These $a$, $b$, $c$, and $d$ are the coefficients of general cubic equation (\ref{cubiceq}). From (\ref{discriminant}), we have
\begin{gather}\label{deltastar}
\Delta = 18 b c d - 4 b^3 d + b^2 c^2 - 4 c^3 - 27 d^2 = d ((18 c - 4 b^2)b - 27 d) + c^2 (b^2 - 4c). 
\end{gather}

From (\ref{stablecond}), we have $c = \beta^2 \epsilon^2 + \beta \mu^2 (\beta \gamma^2 - 2 \epsilon) + \alpha \beta \epsilon \mu^2  > 0$ (note that: $\beta \gamma^2 - 2 \epsilon > 0$) and $d < 0$. Hence, to show that $\Delta < 0$ in equation (\ref{deltastar}), it will suffice to show
\begin{gather}\label{disccond}
\begin{cases}
b > 0 \\
b^2 - 4c < 0. 
\end{cases}
\end{gather}

We will show that condition (\ref{stablecond}) implies (\ref{disccond}). We have
\begin{gather*}
b = (\alpha - 1) \mu^2 - \beta^2 \gamma^2 + 2 \beta \epsilon + \alpha \beta \gamma \mu > (\alpha - 1) \mu^2 - \beta^2 \gamma^2 + \alpha \beta \gamma \mu + \beta^2 \gamma^2 -  \alpha \beta \gamma \mu + (1 - \alpha) \mu^2 = 0 \\
\left[\text{Note that} \ \epsilon > \frac{\beta \gamma^2}{2} - \left(\frac{\alpha \gamma \mu}{2} - \frac{(1 - \alpha) \mu^2}{2 \beta}\right) \right],
\end{gather*}

and
\begin{gather*}
b^2 - 4c = \alpha^2 \beta^2 \gamma^2 \mu^2 + 2 \alpha^2 \beta \gamma \mu^3 + \alpha^2 \mu^4 - 2 \alpha \beta^3 \gamma^3 \mu - 2 \alpha \beta^2 \gamma^2 \mu^2 + \\
+ 4 \epsilon \alpha \beta^2 \gamma \mu - 2 \alpha \beta \gamma \mu^3 - 2 \alpha \mu^4 + \beta^4 \gamma^4 - 4 \epsilon \beta^3 \gamma^2 - 2 \beta^2 \gamma^2 \mu^2 + 4 \epsilon \beta \mu^2 + \mu^4 = \\
= (\alpha - 1)^2 \mu^4 + \beta \mu^2 (\alpha^2 \beta \gamma^2 - 2 \alpha \beta \gamma^2 - 2 \beta \gamma^2 + 4 \epsilon) + 2 \alpha \beta \gamma \mu^3 (\alpha - 1) + \\
+ \alpha \beta^2 \gamma \mu (-2 \beta \gamma^2 + 4 \epsilon) + \beta^3 \gamma^2 (\beta \gamma^2 - 4 \epsilon) \stackrel{\text{(a)}}{<} \\
< (\alpha - 1)^2 \mu^4 + \beta \mu^2 (\alpha^2 \beta \gamma^2 - 2 \alpha \beta \gamma^2 - 2 \beta \gamma^2 + 4 \epsilon) + 2 \alpha \beta \gamma \mu^3 (\alpha - 1) = \\
= (\alpha - 1) \mu^3 ((\alpha - 1) \mu + \alpha \beta \gamma) + \alpha \beta \gamma \mu^3 (\alpha - 1) + \beta \mu^2 (\alpha \beta \gamma^2 (\alpha - 2) - 2 (\beta \gamma^2 - 2 \epsilon)) \stackrel{\text{(b)}}{<} \\
< (\alpha - 1) \mu^3 ((\alpha - 1) \mu + \alpha \beta \gamma). 
\end{gather*}

(where in (a) and (b) we use the facts that $\frac{\beta \gamma^2}{4} < \epsilon < \frac{\beta \gamma^2}{2} \Leftrightarrow \beta \gamma^2 - 4 \epsilon < 0 < \beta \gamma^2 - 2 \epsilon$). \\

From condition (\ref{stablecond}), we have $(\alpha - 1) \mu + \alpha \beta \gamma > (\alpha - 1) \mu + (1 - \alpha) \mu = 0$. Note that $\gamma > \frac{(1 - \alpha) \mu}{\alpha \beta} \Leftrightarrow \alpha \beta \gamma > (1 - \alpha) \mu$. Therefore, we have $b > 0$ and $b^2 - 4c < 0$. Hence, $A^{+} A^{-}$ has no negative real eigenvalues under condition (\ref{stablecond}). \\

By Lemma \ref{lemmaaplus}, $A^{+}$ is Hurwitz for all positive $\beta$, $\gamma$, $\mu$, $\epsilon$ and $\alpha \in (0,1)$. By Lemma \ref{aminuscond}, $A^{-}$ is Hurwitz under condition (\ref{stablecond}). It is easy to verify that the difference $A^{+} - A^{-}$ has rank one. $A^{+} A^{-}$ has no negative real eigenvalues under condition (\ref{stablecond}). Hence, by Proposition \ref{prop4}, $u^\prime(t) = A^{+} u(t)$ and $u^\prime(t) = A^{-} u(t)$ have a CQLF. Therefore, by Proposition \ref{prop1}, the system (\ref{theorem2_system}) is exponentially stable under condition (\ref{stablecond}). This completes the proof. $\Box$ \\

As a useful corollary of Lemma \ref{lemmaaminus}, we have the following fact. 

\begin{cor}\label{lemmaalpha0}
If matrix $A^{-}$ in \emph{(\ref{aminus})} is Hurwitz for some positive $\beta$, $\gamma$, $\mu$, $\epsilon$, and $\alpha \in (0,1)$, then it remains Hurwitz if $\alpha$ is replaced by any $0 < \alpha_0 \leq \alpha$. 
\end{cor}

\textit{Proof}. From (\ref{aminushurwitz}), for any $\alpha_0 \in (0,\alpha]$, we have
\begin{gather}\label{alpha0}
\left(\frac{\beta \gamma}{\mu} + (1 - \alpha_0)\right)\left(\frac{\gamma \mu}{\epsilon} + 1\right) \geq \left(\frac{\beta \gamma}{\mu} + (1 - \alpha)\right)\left(\frac{\gamma \mu}{\epsilon} + 1\right) > 1. 
\end{gather}

Application of Lemma \ref{lemmaaminus} completes the proof. $\Box$

\section{Proof of Theorem \ref{thrm3}}\label{theorem3proof}

We also use the machinery of switched linear systems and common quadratic Lyapunov functions (CQLF) to approach the stability of fluid limits. 

\begin{lem}\label{aminuscond2}
For $\beta > 0$, $\mu > 0$ and $\alpha \in (0,1)$, there exists a pair of $\gamma > 0$ and $\epsilon > 0$ satisfying condition
\begin{gather}\label{lemma6}
\begin{cases}
\epsilon < \frac{\beta \gamma^2}{2} - \frac{\alpha \gamma \mu}{2} \\
\gamma > \frac{\alpha \mu}{\beta}. 
\end{cases}
\end{gather}

Moreover, condition \emph{(\ref{lemma6})} implies matrix $A^{-}$ being Hurwitz. 
\end{lem}

\textit{Proof}. For $\beta > 0$, $\mu > 0$ and $\alpha \in (0,1)$, we have $\frac{\alpha \mu}{\beta} > 0$. Hence, we can always find a value of $\gamma > 0$ satisfying the second condition of (\ref{lemma6}). And from the second condition of (\ref{lemma6}), we have 
\begin{gather}
\frac{\beta \gamma^2}{2} - \frac{\alpha \gamma \mu}{2} > 0. 
\end{gather}

Hence, we can always find a value of $\epsilon > 0$ satisfying the first condition of (\ref{lemma6}). By Lemma \ref{lemmaaminus}, condition (\ref{lemma6}) imply matrix $A^{-}$ being Hurwitz. $\Box$ \\

\textit{Conclusion of the proof of Theorem \ref{thrm3}}. From the help of MATLAB symbolic calculation, we have
\begin{gather}\label{aplusinv}
(A^{+})^{-1} = \left( \begin{array}{ccc}
0 & -\frac{(\alpha - 2)}{2 \beta} & \frac{\alpha}{2 \beta} \\
-\frac{1}{\epsilon} & -\frac{\gamma}{\epsilon} & 0 \\
-\frac{1}{\epsilon} & \frac{(\epsilon - \gamma \mu)}{\epsilon \mu} & -\frac{1}{\mu} \end{array} \right)
\end{gather}

and
\begin{gather}
\label{detaplusinv} \det\left((A^{+})^{-1}\right) = -\frac{1}{\beta \epsilon \mu} < 0.
\end{gather} 

Therefore, matrix $(A^{+})^{-1}$ is non-singular. By Proposition \ref{prop5}, to demonstrate that the product $A^{+} A^{-}$ has no negative eigenvalues under condition (\ref{stablecond2}), it will suffice to show that $[(A^{+})^{-1} + \tau A^{-}]$ is non-singular for all $\tau \geq 0$. We have
\begin{gather}
\det[(A^{+})^{-1} + \tau A^{-}] = \nonumber \\ \label{fractiondetAAA} = -\frac{[\beta^2\epsilon^2\mu^2\tau^3 + (\beta^2\epsilon^2 + \beta^2\gamma^2\mu^2 - 2\beta\epsilon\mu^2 + \alpha\beta\epsilon\mu^2)\tau^2 + (\mu^2 - \alpha\mu^2 + \beta^2\gamma^2 - 2\beta\epsilon - \alpha\beta\gamma\mu)\tau + 1]}{\beta\epsilon\mu}
\end{gather}

(Expression (\ref{fractiondetAAA}) is also obtained with the help of MATLAB symbolic calculation.) To show $\det[(A^{+})^{-1} + \tau A^{-}] \neq 0$ for all $\tau \geq 0$, it will suffice to show numerator of the fraction (\ref{fractiondetAAA}) is not equal to 0. For all $\tau \geq 0$, we have
\begin{gather*}
\beta^2\epsilon^2\mu^2\tau^3 + (\beta^2\epsilon^2 + \beta^2\gamma^2\mu^2 - 2\beta\epsilon\mu^2 + \alpha\beta\epsilon\mu^2)\tau^2 + (\mu^2 - \alpha\mu^2 + \beta^2\gamma^2 - 2\beta\epsilon - \alpha\beta\gamma\mu)\tau + 1 > \\
> (\beta^2\epsilon^2 + (\beta^2\gamma^2 - 2\beta\epsilon)\mu^2 + \alpha\beta\epsilon\mu^2)\tau^2 + ((1 - \alpha)\mu^2 + \beta^2\gamma^2 - 2\beta\epsilon - \alpha\beta\gamma\mu)\tau > 0 \\
\left[\text{Note that} \ \epsilon < \frac{\beta \gamma^2}{2} - \frac{\alpha \gamma \mu}{2} \Leftrightarrow \beta^2\gamma^2 - 2\beta\epsilon - \alpha\beta\gamma\mu > 0 \Rightarrow \beta^2\gamma^2 - 2\beta\epsilon > 0 \right].
\end{gather*}

This implies that $A^{+} A^{-}$ has no negative eigenvalues under condition (\ref{stablecond2}). Hence, $u^\prime(t) = A^{+} u(t)$ and $u^\prime(t) = A^{-} u(t)$ have a CQLF. Therefore, the system (\ref{theorem2_system}) is exponentially stable under condition (\ref{stablecond2}). This completes the proof. $\Box$ 

\section{Numerical examples}\label{numerical}

In this section, we present some numerical examples to show the good performance of the scheme. Later, we also provide some conjectures based on a variety of simulations. 

\begin{exmp}\label{exmp1}
We use the following set of parameters, which satisfies the condition (\ref{stablecond}) but does not satisfy the condition (\ref{stablecond2}): 
\begin{gather*}
\Lambda = 1000 \ , \ \alpha = 0.7 \ , \ \beta = 1 \ , \ \mu = 1 \ , \ \gamma = 2 \ , \ \epsilon = 1.5
\end{gather*}
\end{exmp}

\begin{figure}[h]
    \centering
    \begin{subfigure}[b]{0.45\textwidth}
        \includegraphics[width=\textwidth]{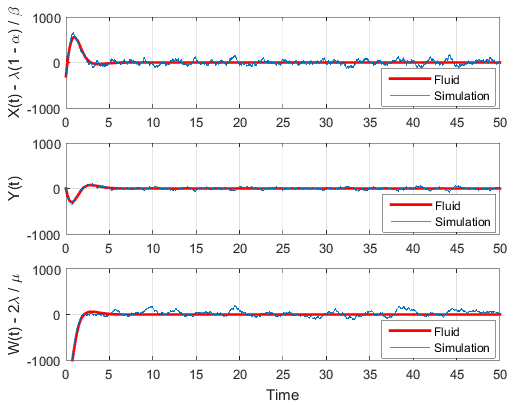}
        \caption{$(X(0),Y(0),Z(0)) = (0,0,0)$}
        \label{exp01a}
    \end{subfigure}
    %~ %add desired spacing between images, e. g. ~, \quad, \qquad, \hfill etc. 
    \begin{subfigure}[b]{0.45\textwidth}
        \includegraphics[width=\textwidth]{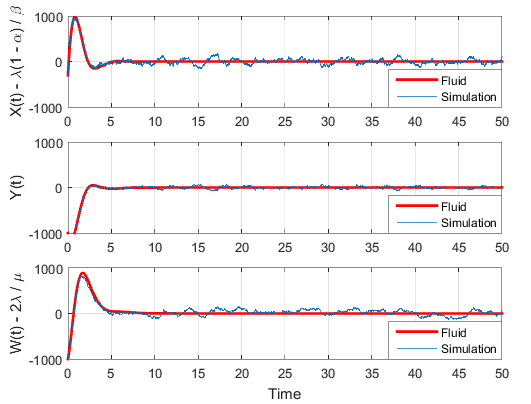}
        \caption{$(X(0),Y(0),Z(0)) = (0,-1000,0)$}
        \label{exp01b}
    \end{subfigure}
    \caption{Comparison of fluid approximations with simulations in Example \ref{exmp1}}\label{exp01}
\end{figure}

We consider two initial conditions: (a) $(X(0),Y(0),Z(0)) = (0,0,0)$; (b) $(X(0),Y(0),Z(0)) = (0,-1000,0)$ (Figure \ref{exp01}). The red line of the figure is the fluid approximation and the blue line of the figure is the simulation experiment. We also did the numerical/simulation experiments with 10 different initial conditions of this set. The results, including those not shown on Figure \ref{exp01}, suggest the global stability of our system. 

\begin{exmp}\label{exmp2}
Let us consider a case when trajectory hits the boundary on $x$. We consider the following set of parameters, which satisfies the condition (\ref{stablecond}) but does not satisfy the condition (\ref{stablecond2}): 
\begin{gather*}
\Lambda = 1000 \ , \ \alpha = 0.5 \ , \ \beta = 3 \ , \ \mu = 2 \ , \ \gamma = 1 \ , \ \epsilon = 1.4
\end{gather*}
\end{exmp}

\begin{figure}[h]
    \centering
    \begin{subfigure}[b]{0.45\textwidth}
        \includegraphics[width=\textwidth]{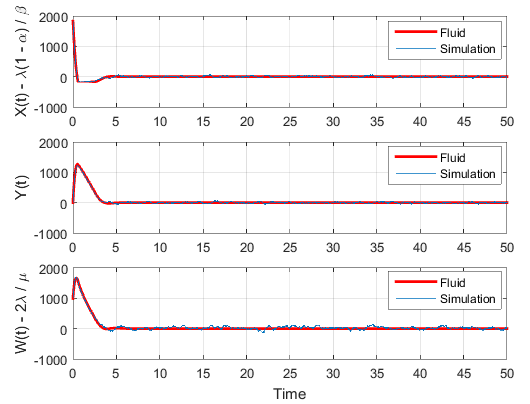}
        \caption{$(X(0),Y(0),Z(0)) = (2000,0,1000)$}
        \label{exp02a}
    \end{subfigure}
    %~ %add desired spacing between images, e. g. ~, \quad, \qquad, \hfill etc. 
    \begin{subfigure}[b]{0.45\textwidth}
        \includegraphics[width=\textwidth]{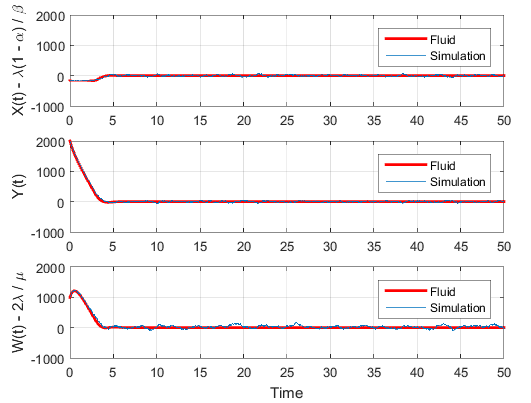}
        \caption{$(X(0),Y(0),Z(0)) = (0,2000,0)$}
        \label{exp02b}
    \end{subfigure}
    \caption{Comparison of fluid approximations with simulations in Example \ref{exmp2}}\label{exp02}
\end{figure}

We consider two initial conditions: (a) $(X(0),Y(0),Z(0)) = (2000,0,1000)$; (b) $(X(0),Y(0),Z(0)) = (0,2000,0)$ (Figure \ref{exp02}). We also did the numerical/simulation experiments with 10 different initial conditions of this set. The results, including those not shown on Figure \ref{exp02}, suggest the global stability of our system even though sometimes the trajectory hits the boundary on $x$. 

\begin{exmp}\label{exmp3}
We use 4 sets of parameters (with different values of $\alpha$), which do not satisfy the condition (\ref{stablecond}) but satisfy the condition (\ref{stablecond2}):  
\begin{gather*}
\Lambda = 1000 \ , \ \beta = 1 \ , \ \mu = 2 \ , \ \gamma = 2 \ , \ \epsilon = 0.19
\end{gather*}
\end{exmp}

\begin{figure}[h]
    \centering
    \begin{subfigure}[b]{0.45\textwidth}
        \includegraphics[width=\textwidth]{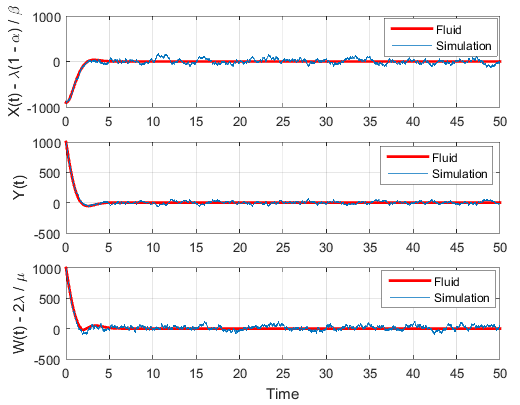}
        \caption{$\alpha_1 = 0.1$}
        \label{exp03a}
    \end{subfigure}
    %~ %add desired spacing between images, e. g. ~, \quad, \qquad, \hfill etc. 
    \begin{subfigure}[b]{0.45\textwidth}
        \includegraphics[width=\textwidth]{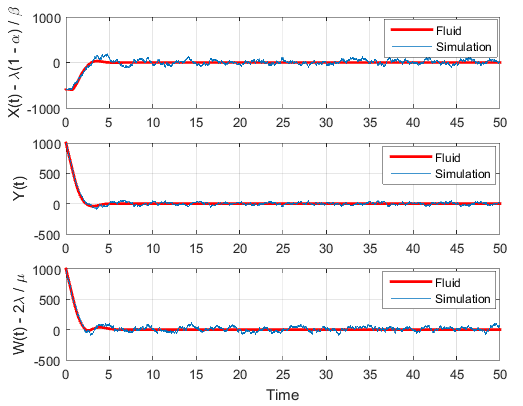}
        \caption{$\alpha_2 = 0.4$}
        \label{exp03b}
    \end{subfigure}
    \begin{subfigure}[b]{0.45\textwidth}
        \includegraphics[width=\textwidth]{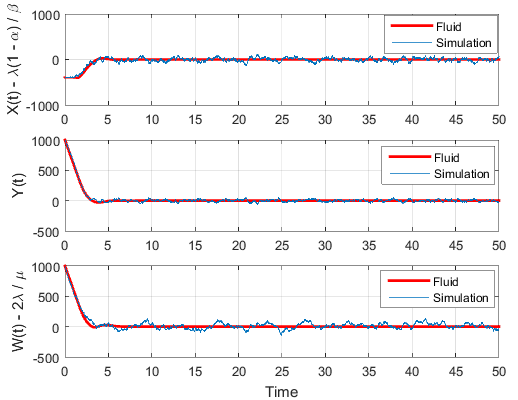}
        \caption{$\alpha_3 = 0.6$}
        \label{exp03c}
    \end{subfigure}
    \begin{subfigure}[b]{0.45\textwidth}
        \includegraphics[width=\textwidth]{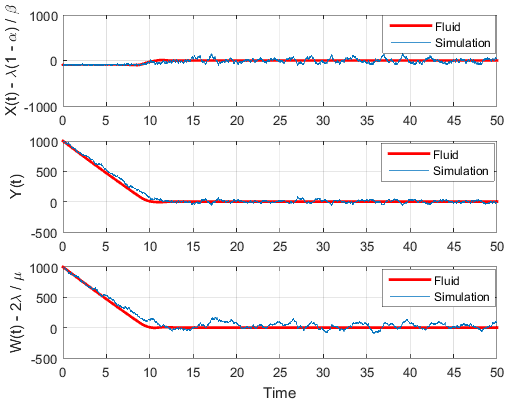}
        \caption{$\alpha_4 = 0.9$}
        \label{exp03d}
    \end{subfigure}
    \caption{Comparison of fluid approximations with simulations in Example \ref{exmp3}}\label{exp03}
\end{figure}

We consider an initial condition $(X(0),Y(0),Z(0)) = (0,1000,500)$ with 4 different values of $\alpha$ ($\alpha_1 = 0.1$, $\alpha_2 = 0.4$, $\alpha_3 = 0.6$, and $\alpha_4 = 0.9$) (Figure \ref{exp03}). We also did the numerical/simulation experiments with 5 different initial conditions for each of the 4 values of $\alpha$. The results, including those not shown on Figure \ref{exp03}, suggest the global stability of our system even though sometimes the trajectory hits the boundary on $x$. \\

Besides these three examples, we also ran the numerical/simulation experiments with another 5 sets of parameters as well as many other different initial conditions of these sets, which satisfy either the condition (\ref{stablecond}) or (\ref{stablecond2}). All these results still suggest the global stability of our system even though sometimes the trajectory hits the boundary on $x$. 

\begin{exmp}\label{exmp4}
In this example, we use a set of parameters, which satisfies neither the condition (\ref{stablecond}) nor (\ref{stablecond2}), but $A^{-}$ is Hurwitz: 
\begin{gather*}
\Lambda = 1000 \ , \ \alpha = 0.5 \ , \ \beta = 1 \ , \ \mu = 2 \ , \ \gamma = 2 \ , \ \epsilon = 3
\end{gather*}
\end{exmp}

\begin{figure}[h]
    \centering
    \begin{subfigure}[b]{0.45\textwidth}
        \includegraphics[width=\textwidth]{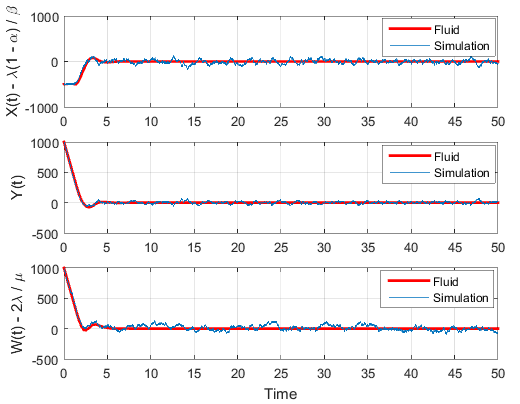}
        \caption{$(X(0),Y(0),Z(0)) = (0,1000,500)$}
        \label{exp04a}
    \end{subfigure}
    %~ %add desired spacing between images, e. g. ~, \quad, \qquad, \hfill etc. 
    \begin{subfigure}[b]{0.45\textwidth}
        \includegraphics[width=\textwidth]{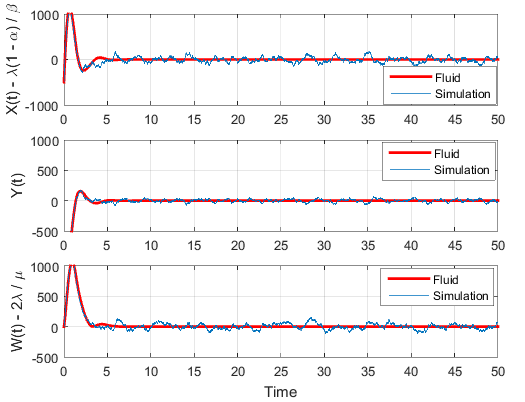}
        \caption{$(X(0),Y(0),Z(0)) = (0,-1000,0)$}
        \label{exp04b}
    \end{subfigure}
    \caption{Comparison of fluid approximations with simulations in Example \ref{exmp4}}\label{exp04}
\end{figure}

We consider two initial conditions: (a) $(X(0),Y(0),Z(0)) = (0,1000,500)$; (b) $(X(0),Y(0),Z(0)) = (0,-1000,0)$ (Figure \ref{exp04}). Besides this example, we also did the numerical/simulation experiments with 5 sets of parameters as well as many other different initial conditions of these sets, which satisfy neither the condition (\ref{stablecond}) nor (\ref{stablecond2}), but $A^{-}$ is Hurwitz. All these results, including those not shown on Figure \ref{exp04}, suggest the local and global stability of our system even though sometimes the trajectory hits the boundary on $x$. 

\begin{exmp}\label{exmp5}
Let us consider the case when $A^{-}$ is not Hurwitz. We use the following two sets of parameters:  
\begin{gather*}
\text{(a)} \ \Lambda = 1000 \ , \ \alpha = 0.5 \ , \ \beta = 0.05 \ , \ \mu = 0.5 \ , \ \gamma = 1 \ , \ \epsilon = 1 \ , \ \text{and} \\
\text{(b)} \ \Lambda = 1000 \ , \ \alpha = 0.9 \ , \ \beta = 0.05 \ , \ \mu = 0.5 \ , \ \gamma = 1 \ , \ \epsilon = 1
\end{gather*}
\end{exmp}

\begin{figure}[h]
    \centering
    \begin{subfigure}[b]{0.45\textwidth}
        \includegraphics[width=\textwidth]{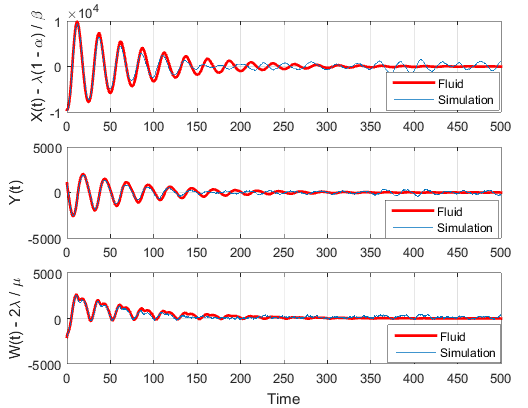}
        \caption{$\alpha = 0.5$}
        \label{exp05a}
    \end{subfigure}
    %~ %add desired spacing between images, e. g. ~, \quad, \qquad, \hfill etc. 
    \begin{subfigure}[b]{0.45\textwidth}
        \includegraphics[width=\textwidth]{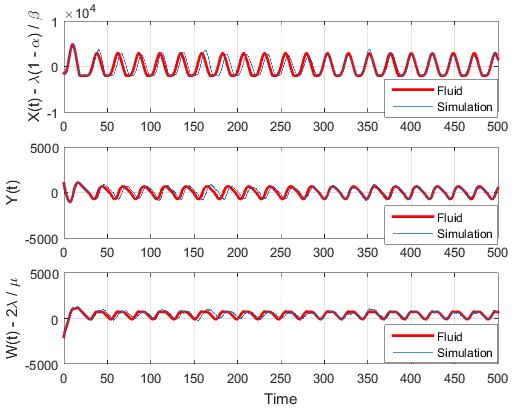}
        \caption{$\alpha = 0.9$}
        \label{exp05b}
    \end{subfigure}
    \caption{Comparison of fluid approximations with simulations in Example \ref{exmp5}}\label{exp05}
\end{figure}

The only difference between these two sets is the parameter $\alpha$. We consider an initial condition $(X(0),Y(0),Z(0)) = (500,1000,500)$ (Figure \ref{exp04}). We see a converging trajectory on the Figure \ref{exp05a} (on the left); in fact, we see convergence for a large number of other initial conditions, for the same set of parameters. Figure \ref{exp05b} shows a trajectory that never converges, under a different set of parameters. \\

The results of Examples \ref{exmp1}, \ref{exmp2} and \ref{exmp3} suggest that our system is globally stable under either the condition (\ref{stablecond}) or (\ref{stablecond2}). The results of Example \ref{exmp4} suggest that our system is locally and globally stable when $A^{-}$ is Hurwitz even if neither the condition (\ref{stablecond}) nor (\ref{stablecond2}) is satisfied. The results of Example \ref{exmp5} suggest that our system might be globally stable under some sets of parameters, but unstable under some different sets of parameters, when $A^{-}$ is not Hurwitz. The summary of our conjectures, motivated by the numerical/simulation experiments, is as follows: 

\begin{conj}
Our system is globally stable if it is locally stable. 
\end{conj}

\begin{conj}
Matrix $A^{-}$ being Hurwitz is sufficient for local stability of our system. ($A^{+}$ is always Hurwitz in our case.)
\end{conj}

\begin{conj} 
If $A^{-}$ is not Hurwitz, the system may be locally stable or locally unstable depending on the parameters.
\end{conj}

\section{Conclusions}\label{conclusion}

In this paper, we study a feedback-based agent invitation scheme for a model with randomly behaving agents. This model is motivated by a variety of existing and emerging applications.
The focus of the paper is on the stability properties of the system fluid limits, arising as asymptotic limits of the system process, when the system scale (customer arrival rate) grows to infinity. The dynamic system, describing the behavior of fluid limit trajectories has a very complex structure -- it is a switched linear system, which in addition has a reflecting boundary. We derived some sufficient local stability conditions, using the machinery of switched linear systems and common quadratic Lyapunov functions. Our simulation and numerical experiments show good overall performance of the feedback scheme, when the local stability conditions hold. They also suggest that, for our model, the local stability is in fact sufficient for the global stability of fluid limits. Verifying these conjectures, as well as expanding the sufficient local stability conditions, is an interesting subject for future research. Further generalizations of the agent invitation model are also of interest from both theoretical and practical points of view.

\bibliographystyle{abbrv}
\bibliography{reference}

\end{document}